\documentclass[12pt]{article}
\usepackage{fancyhdr}
\usepackage{graphicx}
\usepackage{geometry}

\usepackage[T1]{fontenc}
\usepackage[latin1]{inputenc}

\def \bop {\noindent\textbf{Proof. }}
\def \eop {\hbox{}\nobreak\hfill
\vrule width 2mm height 2mm depth 0mm
\par \goodbreak \smallskip}
\newcommand{\1}{1\!\!1}

\def\figurename{Figure} 
\makeatletter
\renewcommand{\fnum@figure}[1]{\figurename~\thefigure.}
\makeatother

\def\tablename{Table} 
\makeatletter
\renewcommand{\fnum@table}[1]{\tablename~\thetable.}
\makeatother
\newcommand{\eps}  {\varepsilon}

\usepackage{amsmath}
\usepackage{amssymb}
\usepackage{amsfonts}
\usepackage{amsthm,amscd}

\newtheorem{theorem}{Theorem}[section]
\newtheorem{lemma}[theorem]{Lemma}
\newtheorem{corollary}[theorem]{Corollary}
\newtheorem{proposition}[theorem]{Proposition}
\theoremstyle{definition}
\newtheorem{definition}[theorem]{Definition}

\theoremstyle{remark}
\newtheorem{remark}[theorem]{Remark}

\numberwithin{equation}{section}
\def\N{\mbox{I\hspace{-.15em}N}}
\def\R{\mbox{I\hspace{-.15em}R}}
\def\P{\mbox{I\hspace{-.15em}P}}
\def\E{\mbox{I\hspace{-.15em}E}}

\begin{document}


\title{ {\bf Homogenization of semilinear PDEs with discontinuous averaged coefficients}}
\author{\quad K. Bahlali$^1$\; A. Elouaflin$^2$\; E. Pardoux$^3$}
\date{}
 \maketitle


\begin{center}
{$^{1}$ IMATH, UFR Sciences, USVT, B.P. 132, 83957 La Garde Cedex, France.
 \\ \small{e-mail:
bahlali@univ-tln.fr}}
\\
{$^{2}$ UFRMI, Universit\'e de Cocody,  22 BP 582 Abidjan, C\^{o}te
d'Ivoire \\
\small{e-mail:  elabouo@yahoo.fr}} \\
{$^3$ LATP, CMI Universit\'e de Provence, 39 rue
Joliot-Curie,13453 Marseille \\
\small{e-mail: pardoux@cmi.univ-mrs.fr } }
\end{center}
\footnotetext[1]{\ Partially supported by PHC Tassili 07MDU705 and  Marie Curie ITN, no. 213841-2.}
\footnotetext[2]{\ Supported by AUF bourse post-doctorale 07-08, R\'ef.: PC-420/2460.}


\pagestyle{fancy} \fancyhead{} \fancyhead[EC]{K. Bahlali, A.
Elouaflin and E. Pardoux} \fancyhead[EL,OR]{\thepage}
\fancyhead[OC]{Homogenization of semilinear PDEs with discontinuous
coefficients} \fancyfoot{}
\renewcommand\headrulewidth{0.5pt}

\begin{abstract}
We study the asymptotic behavior of solutions of semilinear PDEs. Neither periodicity nor ergodicity will be assumed. On the other hand, we assume that the
 coefficients  have averages in the Cesaro
sense. In such a case, the averaged coefficients could be
discontinuous. We use a probabilistic approach based  on weak
convergence of the
associated backward stochastic differential equation (BSDE) in the Jakubowski S-topology to derive the averaged PDE. However, since the averaged coefficients are discontinuous, the classical viscosity solution is not defined for the averaged PDE. We then use the notion of "$L^p-$viscosity solution"  introduced in \cite{CCKS}. The existence of $L^p-$viscosity solution to  the averaged PDE is proved here by using BSDEs techniques.
\end{abstract}

\noindent {\bf Keys words}: \textit{Backward stochastic differential equations (BSDEs),
$L^p$-viscosity solution for PDEs, homogenization, Jakubowski S-topology, limit in the Cesaro
sense.}

\section{Introduction}
Homogenization of a partial differential equation (PDE) is the
process of replacing rapidly varying coefficients by new ones such
that the solutions are close.
 Let for example $a$ be a one dimensional periodic
  function which is positive and bounded away from zero. For $\varepsilon>0$, we consider the operator
  $$
  L_\varepsilon = div(a(\frac{x}{\varepsilon})\nabla)
  $$
  For small $\varepsilon$, $L_\varepsilon$ can be replaced by
  $$
  L = div(\overline a\nabla)
  $$
  where $\overline a$ is the averaged (or limit, or effective) coefficient
  associated to $a$. As $\varepsilon$ is small, the solution of the
  parabolic equation
$$
  \partial_t u = L_\varepsilon u, \qquad  u(0,x) = f(x)
  $$
  is close to the corresponding solution with $L_\varepsilon$
  replaced by $L$.

The probabilistic approach to homogenization is one way to prove such results in the periodic or ergodic case. It is
based on the asymptotic analysis of the diffusion process associated
to the operator $L_\varepsilon$. The averaged coefficient $\overline
a$ is then determined as a certain "mean" of $a$ with respect to the
invariant probability measure of the diffusion process associated to
$L$.

There is a vast literature on the homogenization of PDEs with
periodic coefficients, see for example the monographs \cite{BLP, FR,
PAN} and the references therein.
 There also exists a considerable literature on the study of
asymptotic analysis of stochastic differential equations (SDEs) with
periodic structures and its connection with homogenization of second
order partial differential equations (PDEs). In view of the connection between
BSDEs and semilinear PDEs, this probabilistic tool has been used in order to prove homogenization results for certain classes of nonlinear PDEs, see in particular
 \cite{BI, BHP, BH, D, EO, I,
L, EP2, P} and the references therein.
The two classical situations which have been mainly studied are the cases of
deterministic periodic and random stationary coefficients.
This paper is concerned with a different situation, building upon earlier results of
 Khasminskii and Krylov.

In \cite{KK}, Khasminskii \& Krylov consider the averaging of
the following family of diffusions process
\begin{equation}\label{EH0}
 \left\{
 \begin{aligned}
U^{1,\,\varepsilon}_t&=\frac{x_1}{\eps}+\frac{1}{\varepsilon}\int_0^t\varphi(U^{1,\,\varepsilon}_s,\,U_s^{2,\,\varepsilon})dW_s,\\
U^{2,\,\varepsilon}_t&=x_2+\int_0^t b^{(1)}
(U^{1,\,\varepsilon}_s,\,U^{2,\,\varepsilon}_s)ds
+\int_0^t\sigma^{(1)}(U^{1,\,\varepsilon}_s,\,U_s^{2,\,\varepsilon})d\widetilde{W}_s,
\end{aligned}
\right.
\end{equation}
where for each $\varepsilon>0$ small, $U^{1,\,\varepsilon}_t$ is a one-dimensional
null-recurrent fast component and
$U^{2,\,\varepsilon}_t$  is a $d$--dimensional slow component. The function $\varphi$
(resp. $\sigma^{(1)}$, resp. $b^{(1)}$) is $\R$-valued (resp. $\R^{d\times (k-1)}$-valued, resp. $\R^{d}$-valued). $(W,\widetilde{W})$ is a $k$-dimensional standard Brownian
motion whose component $W$ (resp. $\widetilde{W}$) is one dimensional
(resp. ($k$-$1$)-dimensional). Define now $(X^{1,\eps},X^{2,\eps})=
(\eps U^{1,\eps},U^{2,\eps})$. The process $\{X^\eps_t := (X^{1,\eps}_t,X^{2,\eps}_t),\ t\ge0\}$
solves the SDE
\begin{equation}\label{EH01}
 \left\{
\begin{aligned}
X^{1,\,\varepsilon}_t&=x_1+\int_0^t\varphi\left(\frac{X^{1,\,\varepsilon}_s}{\varepsilon},\,X_s^{2,\,\varepsilon}\right)dW_s,\\
X^{2,\,\varepsilon}_t&=x_2+\int_0^t b^{(1)}
\left(\frac{X^{1,\,\varepsilon}_s}{\varepsilon},\,X^{2,\,\varepsilon}_s\right)ds
+\int_0^t\sigma^{(1)}\left(\frac{X^{1,\,\varepsilon}_s}{\varepsilon},\,X_s^{2,\,\varepsilon}\right)d\widetilde{W}_s,
\end{aligned}
\right.
\end{equation}
They define the averaged coefficients as limits in the
Cesaro sense. With the additional assumption that the presumed
SDE limit is weakly unique, they prove that the process
$(X^{1,\,\varepsilon}_t,\,X^{2,\,\varepsilon}_t)$
converges in distribution towards a Markov  diffusion
$(X^1_t,\,X^2_t)$. As a byproduct, they derive the limit behavior
of the linear PDE associated to
$(X^{1,\,\varepsilon}_t,\,X^{2,\,\varepsilon}_t)$, in the case where weak uniqueness  of the limiting PDE holds in the Sobolev space $ W_{d+1,\text{loc}}^{1,2}(\R_+\times\R^d)$  of all funcions $u(t,x)$ defined on $\R_+\times\R^d$ such that both $u$ and all the generalized derivatives $D_t u$, $D_x u$, and $D^2_{xx} u$ belong to $L^{d+1}_{loc} (\R_+\times\R^d)$.

\par In the present note, we extend the results of \cite{KK} to parabolic semilinear PDEs.
 Note that the limiting coefficients can be discontinuous. More
  precisely, we consider the following sequence of semi-linear PDEs, indexed by
$\varepsilon
>0$,
\begin{equation}\label{E1}
\left \{
\begin{aligned}
 \frac{\partial
v^{\varepsilon}}{\partial t}(t,\,x_1,\,x_2)&=
(\mathcal{L}^{\varepsilon} v^{\varepsilon})(t,x_1,\,x_2)
+f(\frac{x_1}{\varepsilon},\,x_2,\,v^{\varepsilon}(t,\,x_1,\,x_2)), \quad t>0\\
v^{\varepsilon}(0,\,x_1,\,x_2)&= H(x_1,\,x_2);\  (x^1,x^2) \in \R\times\R^{d}.
\end{aligned}
\right.
\end{equation}
\begin{equation*}
\mathcal{L}^{\varepsilon}:=a_{00}(\frac{x_1}{\varepsilon},\,x_2)\frac{\partial^2
}{\partial^2 x_1}+\sum_{i,\,j =1}^d
a_{ij}(\frac{x_1}{\varepsilon},\,x_2)\frac{\partial^2 }{\partial
x_{2i}\partial x_{2j}} +\sum_{i=1}^d
b_i^{(1)}(\frac{x_1}{\varepsilon},\,x_2)\frac{\partial }{\partial
x_{2i}},
\end{equation*}
where $\varphi$, $\sigma^{(1)}$ and $b^{(1)}$ are those defined above in equation (\ref{EH0}),
$$
a_{00}:=\frac{1}{2}\varphi^2,
\ \ \ \ \ \ \ a_{ij}:=\frac{1}{2}
(\sigma^{(1)}\sigma^{(1)\,*})_{ij},\,i,\,j=1,\,...,\,d,
$$
and the  real valued  measurable functions $f$ and $H$ are defined
on $\R^{d+1}\times\R$ and $\R^{d+1}$ respectively.

We put
   $$b:=\begin{pmatrix}0\\ b^{(1)}\end{pmatrix}, \ \ \ \ a(x):= \frac{1}{2}(\sigma\sigma^\ast)(x),\quad \text{ with }
 \sigma := \begin{pmatrix}
\varphi& 0\\
0 & \sigma^{(1)}
\end{pmatrix}
.$$

We write
$$B:=\begin{pmatrix}W\\ \widetilde{W}\end{pmatrix} \ \ \ \  \text{ and } \ \ \ \ X^{\varepsilon}:=\begin{pmatrix}X^{1,\,\varepsilon}\\ X^{2,\varepsilon}\end{pmatrix}
  .$$
The  PDE (\ref{E1}) is then connected to the system of SDE -- BSDE
\begin{equation}\label{E2}
\left\{
\begin{aligned}
X^{\varepsilon}_s&=
x+\int_0^s b(\frac{X^{1,\,\varepsilon}_r}{\varepsilon},\,X^{2,\,\varepsilon}_r)dr
+\int_0^s \sigma (\frac{X^{1,\,\varepsilon}_r}{\varepsilon},\,X^{2,\,\varepsilon}_r)dB_r,\\
Y^{\varepsilon}_s&=H(X^{\varepsilon}_t)+\int_s^t
f(\frac{X^{1,\,\varepsilon}_r}{\varepsilon},\,X^{2,\,\varepsilon}_r,\,Y^{\varepsilon}_r)dr-\int_s^t
Z^{\varepsilon}_r\,dM_r^{X^{\varepsilon}},\,\forall\,s\in [0,\,t]
\end{aligned}
\right.
\end{equation}
where $M^{X^{\varepsilon}}$ is the martingale part of the
process $X^{\varepsilon}$ i. e.
$$M^{X^{\varepsilon}}_s=\int_0^s \sigma (\frac{X^{1,\,\varepsilon}_r}{\varepsilon},\,X^{2,\,\varepsilon}_r)dB_r, \quad 0\le s\le t.$$
Note that $Y^\varepsilon_0$ does depend upon the pair $(t,x)$ where
$x$ is the initial condition of the forward SDE part of \eqref{E2},
and $t$ is the final time of the BSDE part of \eqref{E2}. It follows from e. g.
Remark 2.6 in \cite{EP1} that under suitable conditions upon the coefficients
$\{v^{\varepsilon}(t,\,x):=Y^{\varepsilon}_0,\, t\ge0,\ x=(x_1,x_2)\in\R^{d+1}\}$
solves the PDE \eqref{E1}.

 \vskip 0.2cm The aim of the present paper is
\begin{enumerate}
\item  to show that for each $t>0$, $x\in\R^{d+1}$, the sequence of processes
 $(X_s^\varepsilon , Y_s^\varepsilon, \int_s^t
Z^{\varepsilon}_r\,dM_r^{X^{\varepsilon}})_{0\leq s\leq t}$ converges in law to the
process $(X_s,Y_s, \int_s^t Z_r\,dM_r^X)_{0\leq s\leq t}$ which is the unique
solution to the system of SDE -- BSDE
\begin{equation}\label{fbsdebar}
\left\{
\begin{aligned}
X_s&=x+\int_0^s \bar{b}(X_r)dr+\int_0^s\bar{\sigma}(X_r)dB_r,\,0\leq s\leq t.\\
{Y}_s&=H(X_t)+\int_s^t \bar{f}(X_r,\,{Y}_r)dr-\int_s^t
{Z}_rdM^{X}_r, 0\leq s\leq t,
\end{aligned}
\right.
\end{equation}
where $M^X$ is the martingale part of X and $\bar \sigma$, $\bar b$ and $\bar f$ are
respectively the average of $\sigma$, $ b$ and $ f$, in a sense which will be made precise below;
\item deduce from the first result that for each $(t,x)$,
$v^{\varepsilon}(t,x_1,x_2)\longrightarrow v(t,x_1,x_2)$, where $v$ solves the following
averaged PDE in the $L^p$-viscosity sense
\begin{equation}\label{E3}
\left\{
\begin{aligned} \frac{\partial{ v}}{\partial
t}(t,\,x_1,\,x_2)&=(\bar{L}v)(t,x_1,\,x_2)
+\bar{f}(x_1,\,x_2,\,v(t,\,x_1,\,x_2)) \,\quad  t>0,
\\
v(0,x_1,\,x_2)&= H(x_1,\,x_2),
\end{aligned}
\right.
\end{equation}
with $$
\bar{L} = \sum_{i,\,j}\bar{a}_{ij}(x_1,\,x_2)\frac{\partial^2}{\partial
x_i\partial
x_j}+\sum_i\bar{b}_i(x_1,\,x_2)\frac{\partial}{\partial x_i}$$
the averaged operator.
\end{enumerate}

The method used to derive the averaged BSDE is
based on weak convergence in
the ${\bf S}$-topology and is close to that used in \cite{EP2} and \cite{P}. In our
framework, we show that the limiting system of SDE -- BSDE (\ref{fbsdebar}) has a
unique solution. However, due  to the discontinuity of the
coefficients, the classical viscosity solution is not defined for the averaged PDE (\ref{E3}). We then use the notion of "$L^p-$viscosity solution". We use BSDE techniques to establish the existence of $L^p-$viscosity solution for the averaged PDE. The notion of  $L^p$-viscosity solution has been introduced by Caffarelli {\it et al.} in \cite{CCKS} to study fully nonlinear PDEs with measurable coefficients. Note however that although the
notion of a $L^p$-viscosity solution is available for PDEs with
merely measurable coefficients, continuity of the solution is required. In our situation, the lack of $L^2$-continuity property for
the flow $X^x :=(X^{1,\,x},\,X^{2,\,x})$ transfers the difficulty to
the backward one and hence we cannot prove the $L^2$-continuity of
the process $Y$. To overcome this difficulty, we establish weak
continuity for the flow $x\mapsto (X^{1,\,x},\,X^{2,\,x})$ and using
the fact that $Y_0^x$ is deterministic, we derive the continuity
property for $Y_0^x$.

The paper is organized as follows: In section $2$, we make precise some
notations and formulate our assumptions. Our main results are stated in
section $3$.
Section $4$ and $5$ are devoted to the proofs.

\section{Notations and assumptions}
\subsection{Notations}

For a given function $g(x_1,\,x_2)$, we define
 \begin{align*}
 g^{+}(x_2)&:=\lim_{x_1\rightarrow
+\infty}\frac{1}{x_1}\int_0^{x_1}g(t,\,x_2)dt\\
g^{-}(x_2)&:=\lim_{x_1\rightarrow
-\infty}\frac{1}{x_1}\int_0^{x_1}g(t,\,x_2)dt
\end{align*}
The average, in Cesaro sense, of $g$ is defined by
$$g^{\pm}(x_1,\,x_2):=g^+(x_2)1_{\{x_1>0\}}+g^-(x_2)1_{\{x_1\leq
0\}}$$

Let $\rho(x_1,\,x_2):=a_{00}(x_1,\,x_2)^{-1} (=[\frac{1}{2}\varphi^2(x_1,\,x_2)]^{-1})$
and denote by $\bar
{b}(x_1,\,x_2),\,\bar{a}(x_1,\,x_2)$ and $\bar
{f}(x_1,\,x_2,\,y)$, the averaged coefficients defined by
\begin{align*}
\bar
{b}_i(x_1,\,x_2)&=\frac{(\rho b_i)^{\pm}(x_1,\,x_2)}{\rho^{\pm}(x_1,\,x_2)},\qquad i=1,\,...,\,d \\
\bar
{a}_{ij}(x_1,\,x_2)&=\frac{(\rho a_{ij})^{\pm}(x_1,\,x_2)}{\rho^{\pm}(x_1,\,x_2)}, \qquad
\,i,\,j=0,\,1,\,...,\,d\\
\bar
{f}(x_1,\,x_2,\,y)&=\frac{(\rho f)^{\pm}(x_1,\,x_2,\,y)}{\rho^{\pm}(x_1,\,x_2)}. \\
\bar
{\sigma}(x_1,\,x_2)& = (\bar a(x_1,x_2))^\frac{1}{2}
\end{align*}
where $\bar a(x_1,x_2)$ denotes the matrix  $(\bar
{a}_{ij}(x_1,\,x_2))_{i,j}$.

It is worth noting that $\bar{b},\,\bar{a}$ and
$\bar{f}$ may be discontinuous at $x_1=0$.

\subsection{Assumptions.}
We consider the following conditions.
\begin{trivlist}
\item {\bf(A1)} The functions $b^{(1)},\,\sigma^{(1)},\,\varphi$ are uniformly  Lipschitz in the variables
$(x_1,\,x_2)$.
\item {\bf(A2)}
For each $x_1$, the first and second order derivatives
with respect to $x_2$ of these functions
are bounded continuous functions of $x_2$.
\item {\bf(A3)} $a^{(1)}:=\frac{1}{2} (\sigma^{(1)}\sigma^{(1)\,*})$  is uniformly elliptic,
i. e.
$\exists \Lambda >0; \quad \forall x, \xi \in \R^d, \quad  \xi^*a^{(1)}(x)\xi \geq\Lambda \vert \xi \vert^2 $.
Moreover, there exist  positive constants $C_1,\,C_2,\,C_3$ such that
\begin{equation*}
\left\{
\begin{aligned}
( i )&\quad C_1\leq a_{00}(x_1,\,x_2)\leq C_2       \\
(ii )&\quad \vert a^{(1)}(x_1,\,x_2)\vert + \vert b(x_1,\,x_2)\vert^2\leq
C_3(1+|x_2|^2).
\end{aligned}
\right.
\end{equation*}
\item {\bf(B1)} Let
$D_{x_2}\rho$ and $D_{x_2}^2\rho$ denote respectively the gradient vector
and the matrix of second derivatives of $\rho$ with respect to $x_2$.
We assume that uniformly with respect to $x_2$
\begin{align*}
   \frac{1}{x_1}\int_0^{x_1}\rho(t,\,x_2)dt&\longrightarrow
   \rho^{\pm}(x_2)
    \qquad \hbox{as} \quad {x_1\rightarrow\pm\infty},\\
    \frac{1}{x_1}\int_0^{x_1}D_{x_2}\rho(t,\,x_2)dt&\longrightarrow
   D_{x_2}\rho^{\pm}(x_2)
    \qquad \hbox{as} \quad {x_1\rightarrow\pm\infty},\\
       \frac{1}{x_1}\int_0^{x_1}D_{x_2}^2\rho(t,\,x_2)dt&\longrightarrow
   D_{x_2}^2\rho^{\pm}(x_2)
   \qquad \hbox{as} \quad {x_1\rightarrow\pm\infty}.
   \end{align*}
\end{trivlist}

\begin{trivlist}
\item {\bf(B2)}  For every $i$ and $j$, the coefficients $\rho b_i,\,D_{x_2}(\rho b_i),\,D_{x_2}^2(\rho b_i),\,\rho a_{ij},\,D_{x_2}(\rho a_{ij}),\\D_{x_2}^2(\rho a_{ij})$
have  averages in the Cesaro sense.\\
\item {\bf(B3)} For every function
$k\in\{\rho b_i,\,D_{x_2}(\rho b_i),\,D_{x_2}^2(\rho b_i),\,\rho a_{ij},\,D_{x_2}(\rho a_{ij}),\,D_{x_2}^2(\rho a_{ij})\}$, there\\ exists a bounded function $\alpha:\R^{d+1}\to\R$ such that
\begin{equation}\label{G1}
\left\{\begin{aligned}
\frac{1}{x_1}\int_0^{x_1}k(t,\,x_2)dt-k^{\pm}(x_1,\,x_2)&=(1+|x_2|^2)\alpha
(x_1,\,x_2),\\
 \lim_{|x_1|\longrightarrow \infty}\sup_{x_2\in
{\R}^d}|\alpha (x_1,\,x_2)|&=0.
\end{aligned}
\right.
\end{equation}
\end{trivlist}
\begin{trivlist}
\item {\bf(C1)}
\begin{description}
\item{(i)}  The coefficient $f$ is uniformly Lipschitz in $(x_1,x_2,y)$ and, for each $x_1\in\R$, its derivatives in $(x_2,y)$ up to and including second order derivatives are bounded continuous functions of $(x_2,y)$.
\item{(ii)} There exists  positive constant $K$ such that

\hskip 1.0cm  for every $(x_1,x_2,y)$, \qquad
$
 \vert f(x_1,x_2, y)\vert \leq
K(1+|x_2|+\vert y\vert).
$
\item{(iii)} $H$ is continuous and bounded.
\end{description}
\end{trivlist}
\begin{trivlist}
\item {\bf(C2)}   $\rho f$ has a limit in the Cesaro
sense and there exists a bounded measurable function $\beta:\R^{d+2}\to\R$ such that
\begin{equation}\label{G2}
\left\{
\begin{aligned}
\frac{1}{x_1}\int_0^{x_1}\rho(t,\,x_2)f(t,\,x_2,\,y)dt-(\rho f)^{\pm}(x_1,\,x_2,\,y)&=(1+|x_2|^2+|y|^2)\beta
(x_1,\,x_2,\,y)\\
\lim_{|x_1|\rightarrow \infty}\sup_{(x_2,\,y) \in
{\R}^d\times \R}|\beta
(x_1,\,x_2,\,y)|&=0,
\end{aligned}
\right.
\end{equation}

\item {\bf(C3)}\,\,For each $x_1$, $\rho f$ has
derivatives up to second order in $(x_2,y)$
and these derivatives  are bounded  and  satisfy (C2).
\end{trivlist}
\noindent Throughout the paper, {\bf (A)} stands for conditions (A1), (A2), (A3);
{\bf (B)} for conditions (B1), (B2), (B3) and {\bf (C)} for  (C1), (C2), (C3).\\

\section{The main results}

Consider the equation
\begin{equation}\label{EH1}
X_t^{x}=x+\int_0^t \bar{b}(X_s^{x})ds+
\int_0^t \bar{\sigma} (X_s^{x})dB_s,\,t\ge0.
\end{equation}
Assume that {\bf (A), (B)} hold. Then, from Khasminskii \& Krylov \cite{KK} and Krylov \cite{K}, we deduce that for each fixed, $x\in\R^{d+1}$  {\it the process $X^{\varepsilon}:=(X^{1,\,\varepsilon},\,X^{2,\,\varepsilon})$
converges in distribution to the process $X:=(X^1,\,X^2)$} which is the unique weak solution to SDE (\ref{EH1}).

\vskip 0.2cm We now define the notion of $L^p$-viscosity solution of a parabolic PDE.  This notion has been introduced by Caffarelli {\it et al.}
 in \cite{CCKS} to study PDEs with measurable coefficients. Presentations of this topic can be found in
\cite{CCKS, CKLS}.

\noindent Let $g\,: \R^{d+1}\times\R\longmapsto\R$
be a measurable function  and
$$
\bar{L}:=\sum_{i,\,j}\bar{a}_{ij}(x)\frac{\partial^2}{\partial
x_i\partial
x_j}+\sum_i\bar{b}_i(x)\frac{\partial}{\partial x_i} $$
denote the second order PDE operator associated to the SDE (\ref{EH1}).

\noindent We consider the parabolic equation
\begin{equation}\label{EH2}
\left\{
\begin{aligned}
\frac{\partial{ v}}{\partial
t}(t,x)&=(\bar {L}v)(t,x)
+g(x,v(t,x)),\,\,t\ge0 \\
v(0,x)&= H(x).
\end{aligned}
\right.
\end{equation}

\vskip 0.2cm
\begin{definition}
\label{DED1} Let $p$ be an integer such that $p>d+2$.

(a) \ A function $v\in \mathcal{C}\left([0,\,T]\times
\R^{d+1},\,\R\right)$ is a $L^p$-viscosity sub-solution of the PDE
(\ref{EH2}), if for every $x\in\R^{d+1}$, $v(0,\,x)\leq H(x)$ and
for every $\varphi \in W^{1,\,2}_{p,\,loc}\left(\R_+\times
\R^{d+1},\,\R\right)$ and $(\widehat{t},\,\widehat{x})\in
(0,\,T]\times \R^{d+1}$ at which $v-\varphi$ has a local maximum,
one has
\begin{equation*}
\text{ess}\!\!\liminf_{(t,\,x)\rightarrow
(\widehat{t},\,\widehat{x})}\left\{\frac{\partial{\varphi}}{\partial
t}(t,x) -(\bar{L}\varphi)(t,x)
- g(x,\,v(t,x))\right\}\leq 0.
\end{equation*}

(b) \  A function $v\in \mathcal{C}\left([0,\,T]\times
\R^{d+1},\,\R\right)$ is a $L^p$-viscosity super-solution of the PDE
(\ref{EH2}), if for every $x\in\R^{d+1}$, $v(0,\,x)\geq H(x)$ and
for every $\varphi \in W^{1,\,2}_{p,\,loc}\left(\R_+\times
\R^{d+1},\,\R\right)$ and $(\widehat{t},\,\widehat{x})\in
(0,\,T]\times \R^{d+1}$ at which $v-\varphi$ has
 a local minimum, one has
\begin{equation*}
\text{ess}\!\!\limsup_{(t,\,x)\rightarrow
(\widehat{t},\,\widehat{x})}\left\{\frac{\partial{
\varphi}}{\partial t}(t,x)
-(\bar{L}\varphi)(t,x)
- g(x,\,v(t,x))\right\}\geq 0 .
\end{equation*}
Here, $\displaystyle
G(t,\,x,\,\varphi(s,\,x))$ is merely assumed to be
measurable upon the variable $x=:(x_1,\,x_2)$.

\vskip 0.2cm (c) \  A function $v\in \mathcal{C}\left([0,\,T]\times
\R^{d+1},\,\R\right)$ is a $L^p$-viscosity solution if it is both a
$L^p$-viscosity sub-solution and super-solution.
\end{definition}

\begin{remark} Condition (a) means that
 for every $\varepsilon>0,\,r>0$, there exists a set $A\subset{ B_r(\widehat{t},\,\widehat{x})}$ of
positive measure such that, for every $(s,\,x)\in A$,
$$
 \frac{\partial{\varphi}}{\partial
s}(s,\,x)-(\bar{L}\varphi)(t,x)
- g(x,\,v(t,x))\leq \varepsilon.
$$
\end{remark}


The main results are (the ${\bf S}$--topology is explained in the Appendix below)

\begin{theorem}\label{th1} Assume (A), (B), (C) hold. Then, for any $(t,x)\in\R_+\times\R^{d+1}$, there exists a process $(X_s,Y_s, Z_s)_{0\le s\le t}$ such that,


 \noindent (i) the sequence of process $X^{\varepsilon}$ converges in law to the continuous process X, which is the unique weak solution to SDE (\ref{fbsdebar}), in $C([0,t];\R^{d+1})$ equipped with the uniform topology.

 \noindent (ii) the sequence of processes $(Y_s^\varepsilon, \int_s^t
Z^{\varepsilon}_r\,dM_r^{X^{\varepsilon}})_{0\le s\le t}$ converges in law to the process \\
\
$(Y_s, \int_s^t
Z_r\,dM_r^{X})_{0\le s\le t}$ in $D([0,t];\R^2)$,\ where $M^X$ is the martingale part of $X$,\  equipped with the ${\bf S}$--topology.

\noindent (iii) (Y,Z) is the unique solution to BSDE (\ref{fbsdebar}) such that,

(a) \ \ (Y,Z) is $\mathcal{F}^X-$adapted and $ (Y_s, \int_s^t Z_r\,dM_r^X)_{0\le s\le t}$ is continuous.

(b) \ \ $\E\big(\sup_{0\leq s\leq t}\vert Y_s\vert^2 + \int_0^t \vert Z_r\sigma(X_r)\vert^2 dr \big)<\infty$
\end{theorem}
\noindent The uniqueness means that, if $(Y^{1}, Z^{1})$ and $(Y^{2}, Z^{2})$ are two solutions of BSDE (\ref{fbsdebar}) satisfying (iii) (a)-(b) then,  \
$\E\left(\sup_{0\leq s\leq
t}\left|Y^{1}_s-Y^{2}_s\right|^2
+\int_0^t\left|Z_r^{1}\sigma(X_r)-Z_r^{2}\sigma(X_r)\right|^2 dr\right) = 0$, i. e. since
$\sigma\sigma^\ast$ is elliptic (see {\bf (A3)}), $Y^1_s=Y^2_s$ \ $\forall 0\le s\le t$,
$\P$ a. s., and $Z^1_s=Z^2_s$ $ds\times d\P$ a. e.


\begin{theorem}\label{th2} Assume (A), (B), (C) hold. For $\varepsilon >0$, let $v^\varepsilon$ be the unique solution to the problem (\ref{E1}). Let $(Y_s^{(t,x)})_s$ be the unique solution of the BSDE (\ref{fbsdebar}).
Then
\begin{trivlist}
\item{(i)} Equation (\ref{E3}) has a unique $L^p$-viscosity solution $v$ such that $v(t,x) =Y_0^{(t,x)}$.
\item{(ii)} For every $(t,x)\in\R_+\times\R^{d+1}$, $v^\varepsilon(t,x)\to v(t,x)$, as
$\varepsilon\to0$.
\end{trivlist}
\end{theorem}


\section{Proof of Theorem \ref{th1}.}

In all of this section, $(t,x)\in\R_+\times\R^{d+1}$ is arbitrarily fixed with $t>0$.

Assertion (i) follows from \cite{KK} and \cite{K}. Assertion (iii) can be established as in \cite{EP2, P}.  We shall prove (ii). We first deduce from our assumptions (see in particular {\bf(A3)} which says that the coefficients of the forward SDE part of \eqref{E2} are bounded with respect to their first variable, and grow at most linearly in their second variable)
\begin{lemma}\label{estX}
For all $p\ge1$, there exists  constant $C_p$ such that for all $\eps>0$,
$$\E\left(\sup_{0\le s\le t}[|X^{1,\eps}_s|^p+|X^{2,\eps}_s|^p]\right)\le C_p.$$
\end{lemma}

\subsection{Tightness and convergence for the BSDE.}
\begin{proposition}
\label{R5} There exists a positive constant $C$ such that for all
 $\varepsilon>0$
\begin{equation*}
\E\left(\sup_{0\leq s\leq
t}\left|Y^{\varepsilon}_s\right|^2
+\int_0^t\left|Z_r^{\varepsilon}\sigma(X_r^\eps)\right|^2 dr\right)\leq C.
\end{equation*}
\end{proposition}
\bop
We deduce from It\^o's formula (here and below $\bar{X}^{1,\,\varepsilon}_r
=X^{1,\,\varepsilon}_r/\eps$)
\begin{align*}
|Y^{\varepsilon}_s|^2+\int_0^t\left|Z_r^{\varepsilon}
\sigma(X_r^{\varepsilon})\right|^2 dr)&\leq   |H(X^{\varepsilon}_t)|^2+K\int_s^t
|Y^{\varepsilon}_r|^2dr
+\int_s^t|f(\bar{X}^{1,\,\varepsilon}_r,\,X^{2,\,\varepsilon}_r,\,0)|^2 dr\\
&-2\int_s^t \langle Y^{\varepsilon}_r,\,
Z_r^{\varepsilon}dM^{X^{\varepsilon}}_s\rangle.
\end{align*}
It follows from well known results on BSDEs that we can take the expectation in the above identity (see e. g. \cite{EP1}; note that introducing
stopping times as usual and using Fatou's Lemma would yield \eqref{E12} below). We then deduce from Gronwall's lemma that there exists a positive constant $C$ which does
not depend on $\varepsilon$, such that for every $s\in[0,\,t]$,
\begin{equation*}
{\E}\left (|Y^{\varepsilon}_s|^2\right ) \leq C{\E}\left
(|H(X^{\varepsilon}_t)|^2
+\int_0^t|f(\bar{X}^{1,\,\varepsilon}_r,\,X^{2,\,\varepsilon}_r,\,0)|^2dr
\right)
\end{equation*}
and
\begin{equation}
\label{E12}
 {\E}\left
(\int_0^t\left|Z_r^{\varepsilon}\sigma(X_r^\eps)\right|^2 dr\right) \leq  C{\E}\left (|H(X^{\varepsilon}_t)|^2
+\int_0^t|f(\bar{X}^{1,\,\varepsilon}_r,\,X^{2,\,\varepsilon}_r,\,0)|^2dr\right).
\end{equation}
Combining the last two estimates and  the Burkh{\"o}lder-Davis-Gundy
inequality, we get
\begin{equation*}
{\E}\left (\sup_{0\leq s \leq t}|Y^{\varepsilon}_s|^2
+\frac{1}{2}\int_0^t\left|Z_r^{\varepsilon}\sigma(X_r^\eps)\right|^2 dr
\right) \leq  C{\E}\left
(|H(X^{\varepsilon}_t)|^2+\int_0^t|f(\bar{X}^{1,\,\varepsilon}_r,\,X^{2,\,\varepsilon}_r,\,0)|^2dr\right)
\end{equation*}
In view of  condition $ (C1)$ and Lemma \ref{estX}, the proof is
complete.
\eop

We deduce immediately from Proposition \ref{R5}
\begin{corollary}\label{Co1}
$$\sup_{\eps>0}|Y^\eps_0|<\infty.$$
\end{corollary}
\begin{proposition}
\label{R7} For $\varepsilon >0$, let $Y^\varepsilon$ be the process defined by  equation \eqref{E2} and $M^\varepsilon$ be its martingale part. The sequence
$\left(Y^{\varepsilon},\,M^{\varepsilon}\right)_{\varepsilon>0}$ is  tight in the
space
$\mathcal{D}\left([0,\,t],\,\R\right)\times\mathcal{D}\left([0,\,t],\,\R\right)$
 endowed with the $\bf{S}$-topology.
\end{proposition}
\bop Since
$M^{\varepsilon}$ is a martingale, then by \cite{MZ} or \cite{J}, the Meyer-Zheng tightness criteria is
fulfilled whenever
\begin{equation}
\label{ET1}
\sup_{\varepsilon}\left(CV(Y^{\varepsilon})+\E\left(\sup_{0\leq
s\leq
t}|Y^{\varepsilon}_s|+|M^{\varepsilon}_s|\right)\right)<+\infty.
\end{equation}
where the conditional variation $CV$ is
defined in appendix A.
\par\noindent
>From \cite{PV}, the conditional variation $CV(Y^{\varepsilon})$ satisfies
\begin{equation*}
CV(Y^{\varepsilon}) \leq \E\left
(\int_0^t|f(\bar{X}^{1,\,\varepsilon}_s,
\,X^{2,\,\varepsilon}_s,\,Y^{\varepsilon}_s)|ds\right),
\end{equation*}
Now clearly (\ref{ET1}) follows from $({C1})$, Lemma \ref{estX} and Proposition \ref{R5}.
\eop
\begin{proposition}
\label{R8} There exists $\left({Y},\,{M} \right)$ and a countable
subset $\mathsf{D}$ of $ [0,\,t]$ such that along a subsequence
$\varepsilon_n\to0$,
\begin{trivlist}
\item{(i)} $\displaystyle \left({Y}^{\varepsilon_n},{M}^{\varepsilon_n}
\right)\Longrightarrow \left({Y},{M}\right)$ on
$\mathcal{D}\left([0,\,t],\,\R\right)\times\mathcal{D}\left([0,\,t],\,\R\right)$
endowed with the $\bf{S}$--topology.
\item{(ii)}  The finite dimensional distributions of
$\displaystyle \left({Y_s}^{\varepsilon_n},\,{M}_s^{\varepsilon_n}
\right)_{s\in D^c}$ converge to those of $\left({Y}_s,\,{M}_s\right)_{s\in D^c}$.
\item{(iii)}
$(X^{1,\eps_n},X^{2,\eps_n},Y^{\eps_n}) \Longrightarrow(X^1,X^2,{Y})$ , in the
sense of weak convergence in $C([0,t],\R^{d+1})\times
D([0,t],\R)$, equipped with the product of the uniform
convergence and the ${\mathbf S}$ topology.
\end{trivlist}
\end{proposition}
\bop (i) From Proposition \ref{R7}, the family $\left({Y}^{\varepsilon},\,{M}^{\varepsilon}
\right)_{\varepsilon}$ is tight in
$\mathcal{D}\left([0,\,t],\,\R\right)\times\mathcal{D}\left([0,\,t],\,\R\right)$ endowed with the $\bf{S}$-topology. Hence along a subsequence (still denoted by
$\varepsilon$), $\left({Y}^{\varepsilon},\,{M}^{\varepsilon}
\right)_{\varepsilon}$
converges in law on $\mathcal{D}\left([0,\,t],\,\R\right)\times\mathcal{D}\left([0,\,t],\,\R \right)$ towards a c\`ad-l\`ag process $
\left({Y},\,{M}\right)$.

\noindent (ii) follows from Theorem 3.1 in
Jakubowski \cite{J}.

\noindent(iii) According to Theorem \ref{th1} $(i)$, $(X^{1,\eps},X^{2,\eps})\Longrightarrow
(X^1,X^2)$ in $C([0,t],\R^{d+1})$ equipped with the uniform topology.  From assertion (i),
$(Y^{\eps}_\cdot)_{\eps>0}$ is tight in $\mathcal{D}\left([0,\,t],\,\R\right)$ equipped with the $\bf{S}$--topology. Hence the subsequence $\eps_n$ can be chosen in such a way that (iii) holds.
\eop


\subsection{Identification of the limit finite variation process.}

\begin{proposition}
\label{R11} Let $({Y},{M})$ be any limit process as in
Proposition \ref{R8}. Then\\
$(i)$ for every $ s\in[0,\,t] \setminus
\mathsf{D}$,
\begin{equation} \label{E15}
\left\{
\begin{aligned} {Y}_s&=H(X_t)+\int_s^t
\bar{f}(X^1_r,\,\,X^2_r,{Y})dr
 -({M}_t -{M}_s),\\
 \E&\big(\sup_{0\leq s\leq
t}\big[|{Y}_s|^2+|X^1_s|^2+|X^2_s|^2\big]\big)\leq C;
\end{aligned}
\right.
\end{equation}
$(ii)$ ${M}$ is a $\mathcal{F}_s$-martingale, where
$\mathcal{F}_s:=\sigma\left\{X_r,\,Y_r,\,\quad 0\leq r\leq
s\right\}$ augmented with the $\P$-null sets.
\end{proposition}

\noindent To prove this proposition, we need the following lemmas.
\noindent
\begin{lemma}
\label{R2} Assume {\bf (A), (B)}, (C2) and (C3). For $x_2\in\R^d$, $y\in\R$, let
$V^{\varepsilon}(x_1,\,x_2,y)$ denote the solution of the
following equation:
\begin{equation}\label{G3}
\left\{
\begin{aligned}
a_{00}(\frac{x_1}{\varepsilon},\,x_2)D^2_{x_1}V^{\varepsilon}(x_1,\,x_2,y)
&=f(\frac{x_1}{\varepsilon},\,x_2,\,y)-\bar{f}(x_1,\,x_2,\,y), \quad x_1\in\R,\\
V^{\varepsilon}(0,\,x_2,y)=D_{x_1}V^{\varepsilon}(0,\,x_2,y)&=0.
\end{aligned}
\right.
\end{equation}
Then, for some bounded functions $\beta_1$ and $\beta_2$ satisfying (\ref{G2}),

\noindent (i)\quad $
 D_{x_1}V^{\varepsilon}(x_1,\,x_2,y)=x_1(1+|x_2|^2+|y|^2)\beta_1 (\frac{x_1}{\varepsilon},\,x_2,\,y)
,$\\
and the same is true with  $D_{x_1}V^{\varepsilon}$ replaced by
$D_{x_1}D_{x_2}V^{\varepsilon}$ and $D_{x_1}D_yV^{\varepsilon}$;
\vskip 0.25cm\noindent (ii)\quad
$
V^{\varepsilon}(x_1,\,x_2,y)=x_1^2(1+|x_2|^2+|y|^2)\beta_2
(\frac{x_1}{\varepsilon},\,x_2,\,y)
$, \\
and the same is true with $V^{\varepsilon}$ replaced by $D_{x_2}V^{\varepsilon}$, $D_yV^{\varepsilon}$, $D^2_{x_2}V^{\varepsilon}$, $D^2_yV^{\varepsilon}$ and  $D_{x_2}D_yV^{\varepsilon}$.
\end{lemma}
\bop
 We will adapt the idea of \cite{KK} to our situation. For $\eps > 0$ and  $(z,x_2,y)\in\R^{d+2}$
 we set
$$
F_\eps(z,\,x_2,y):=\frac{1}{\eps z}\int_0^{\eps z}\rho(\frac{t}{\varepsilon},\,x_2)g(\frac{t}{\varepsilon},\,x_2,\,y)dt
$$
where
$g(z,\,x_2,\,y):=f(z,\,x_2,\,y)-\bar{f}(\eps z,\,x_2,\,y)$.\\
We only treat the case where $ x_1>0$. The same argument can be used in the case $x_1<0$. We successively use the definition of $\bar f$ and assumptions {\bf(C2)}, to obtain
\begin{align*}
F_\eps(\frac{x_1}{\varepsilon},\,x_2,y)&=\frac{1}{x_1}\int_0^{x_1}\rho(\frac{t}{\varepsilon},\,x_2)f(\frac{t}{\varepsilon},\,x_2,\,y)dt
-(\rho f)^+(x_2,\,y)\\
&\qquad +(\rho f)^+(x_2,\,y)-\frac{(\rho f)^{+}(x_2,\,y)}{\rho^{+}(x_2)}\frac{1}{x_1}\int_0^{x_1}\rho(\frac{t}{\varepsilon},\,x_2)dt\\
&=(1+|x_2|^2+|y|^2)\beta (\frac{x_1}{\varepsilon},\,x_2,\,y)\\
&\quad +\frac{(\rho f)^+(x_2,\,y)}{\rho^+(x_2)}\big[\rho^+(x_2)-\frac{1}{x_1}\int_0^{x_1}
\rho(\frac{t}{\varepsilon},\,x_2)dt\big]\\
&=(1+|x_2|^2+|y|^2)\beta
(\frac{x_1}{\varepsilon},\,x_2,\,y) \\
&\quad +(1+|x_2|^2+|y|^2)\alpha_1(\frac{x_1}{\varepsilon},\,x_2,\,y)
\end{align*}
where $\alpha_1(\frac{x_1}{\varepsilon},\,x_2,\,y) := \frac{(\rho
f)^+(x_2,\,y)}{(1+|x_2|^2+|y|^2)\rho^+(x_2)}\big[\rho^+(x_2)-\frac{1}{x_1}\int_0^{x_1}
\rho(\frac{t}{\varepsilon},\,x_2)dt\big]$. \\
Using assumptions {\bf(B1)} and {\bf(C1-ii)}, one can show that $\alpha_1$ is a  bounded function which satisfies (\ref{G2}).
Since
$D_{x_1}V^{\varepsilon}(x_1,\,x_2,y)=x_1F_\eps(\frac{x_1}{\varepsilon},\,x_2,y)$,
we derive the result for $D_{x_1}V^{\varepsilon}(x_1,\,x_2,y)$.
Further, by integrating it, we get
\begin{align*}
V^{\varepsilon}(x_1,\,x_2,y)&=x_1^2(1+|x_2|^2+|y|^2)\big
((\frac{\varepsilon}{x_1})^2\int_0^{\frac{x_1}{\varepsilon}}t\beta_1 (t,\,x_2,\,y)dt\big),
\end{align*}
where $\beta_1 = \alpha_1 + \beta$.

\noindent Clearly, $ \beta_2
(\frac{x_1}{\varepsilon},\,x_2,\,y):=(\frac{\varepsilon}{x_1})^2
\int_0^{\frac{x_1}{\varepsilon}}t\beta_1 (t,\,x_2,\,y)dt$
 is bounded function which satisfies  (\ref{G2}).  The result for the other quantities can be deduced by
similar arguments from assumptions {\bf(B1)}, {\bf(C1)}, {\bf(C2)} and {\bf(C3)}. \eop

\begin{lemma}\label{R9}
As $\varepsilon \longrightarrow 0$,
$$
\sup_{0\leq s\leq t}\left |\,\int_0^s\left
(f(\frac{X^{1,\,\varepsilon}_r}{\varepsilon},\,X^{2,\,\varepsilon}_r,\,Y^{\varepsilon}_r)
-\bar{f}(X^{1,\,\varepsilon}_r,\,X^{2,\,\varepsilon}_r,\,Y^{\varepsilon}_r)\right
)dr \,\right |\to0$$
in probability  .
\end{lemma}
\bop
We shall show that  for every $s\in [0, \ t]$, $\left
|\int_0^s\big[f(\frac{X^{1,\,\varepsilon}_r}{\varepsilon},\,X^{2,\,\varepsilon}_r,\,
Y^{\varepsilon}_r)-
\bar{f}(X^{1,\,\varepsilon}_r,\,X^{2,\,
\varepsilon}_r,\,Y^{\varepsilon}_r)\big]dr\right
|
$
tends to zero in probability as $\varepsilon$ tends to zero.
Let
$V^{\varepsilon}$
denote the solution of equation (\ref{G3}). Note that $V^{\varepsilon}$ has  first and second
derivatives in $(x_1,\,x_2,y)$ which are possibly discontinuous only at $x_1=0$.
Then, as in \cite{KK},  since $\varphi^2$ is bounded away from zero, we can use
the It{\^o}-Krylov formula to get
\begin{align}\label{ETH1}
&V^{\varepsilon}(X^{1,\,\varepsilon}_s,\,X^{2,\,\varepsilon}_s,Y^\eps_s)
=V^{\varepsilon}(x_1,\,x_2,Y^\eps_0) +\int_0^s
\big[f(\frac{X^{1,\,\varepsilon}_r}{\varepsilon},\,X^{2,\,\varepsilon}_r,\,
Y^{\varepsilon}_r)-
\bar{f}(X^{1,\,\varepsilon}_r,\,X^{2,\,
\varepsilon}_r,\,Y^{\varepsilon}_r)\big]dr \nonumber
\\
&+\int_0^s Trace \big[a^{(1)}(X^{1,\,\varepsilon}_r,\,X^{2,\,\varepsilon}_r)
{D_{x_2}^2V^{\varepsilon}}(X^{1,\, \varepsilon}_r,\,X^{2,\,\varepsilon}_r,Y^\eps_r)\big]dr \nonumber \\
&+\int_0^s
[D_{x_2}V^{\varepsilon}(X^{1,\,\varepsilon}_r,X^{2,\,
\varepsilon}_r,Y^\eps_r)b^{(1)}(X^{1,\,\varepsilon}_r,\,X^{2,\,
\varepsilon}_r) - D_{y} V^\eps(X^{1,\eps}_r,X^{2,\eps}_r,Y^\eps_r)f(\frac{X^{1,\eps}_r}{\eps},
X^{2,\eps}_r,Y^\eps_r)]dr \nonumber \\
&+\int_0^s[D_{x} V^{\varepsilon}(X^{1,\,\varepsilon}_r,\,X^{2,\,\varepsilon}_r,Y^\eps_r)
\sigma (X^{1,\,\varepsilon}_r,\,X^{2,\,\varepsilon}_r)
 +
D_{y} V^\eps(X^{1,\eps}_r,X^{2,\eps}_r,Y^\eps_r)Z^\eps_r
\sigma(\frac{X^{1,\eps}_r}{\eps},X^{2,\eps}_r)]
dB_r \nonumber  \\
&+\frac{1}{2}\int_0^s D_{y}^2 V^\eps(X^{1,\eps}_r,X^{2,\eps}_r,Y^\eps_r)Z^\eps_r
\sigma\sigma^\ast(\frac{X^{1,\eps}_r}{\eps},X^{2,\eps}_r)(Z^\eps_r)^\ast dr \nonumber \\
&+\frac{1}{2}\int_0^s D_{x}D_{y} V^\eps(X^{1,\eps}_r,X^{2,\eps}_r,Y^\eps_r)
\sigma\sigma^\ast(\frac{X^{1,\eps}_r}{\eps},X^{2,\eps}_r)(Z^\eps_r)^\ast dr
\end{align}

In view of Lemma \ref{R2} and Corollary \ref{Co1}, $
V^{\varepsilon}(x_1,\,x_2, Y_0^\eps)$  tends to zero as $\eps\to0$. \\
Using the fact taht $1 = 1_{\{|X^{1,\,\varepsilon}_s| <  \sqrt{\varepsilon}\}}+1_{\{|X^{1,\,\varepsilon}_s|\geq \sqrt{\varepsilon}\}}$ and Lemma \ref{R2}, we obtain
\begin{align*}
\left|V^{\varepsilon}(X^{1,\,\varepsilon}_s,\,X^{2,\,\varepsilon}_s,Y_s^\eps)\right|
&\leq  \varepsilon\left[(1+|X^{2,\,\varepsilon}_s|^2+|Y^{\varepsilon}_s|^2)|
\beta_2(\frac{X^{1,\,\varepsilon}_s}
{\varepsilon},\,X^{2,\,\varepsilon}_s,\,Y^{\varepsilon}_s)|\right]\\
&\  +1_{\{|X^{1,\,\varepsilon}_s|\geq \sqrt{\varepsilon}\}}|X^{1,\,\varepsilon}_s|^2\left[(1+|X^{2,\,\varepsilon}_s|^2
+|Y^{\varepsilon}_s|^2)|\beta_2(\frac{X^{1,\,\varepsilon}_s}{\varepsilon},\,X^{2,\,\varepsilon}_s,\,Y^{\varepsilon}_s)|\right]\\
\end{align*}
{}From Lemma \ref{estX} and Proposition \ref{R5}, we deduce
that
$$
\E\left(\sup_{0\leq s\leq
t}|V^{\varepsilon}(X^{1,\,\varepsilon}_s,\,X^{2,\,\varepsilon}_s, Y_s^\eps)|\right)
\leq
K\left(\varepsilon+\sup_{|x_1|\geq\sqrt{\varepsilon}}
\sup_{(x_2,\,y)}|\beta_2(\frac{x^1}{\varepsilon},\,x^{2},\,y)|\right)
$$
Then, since  $\beta_2$ satisfy respectively
(\ref{G2}), the right hand side of the previous inequality tends to zero as
$\varepsilon\longrightarrow 0$. Similarly, one can show that
\begin{align*}
&\int_0^s Trace \big[a^{(1)}(X^{1,\,\varepsilon}_r,\,X^{2,\,\varepsilon}_r)
{D_{x_2}^2V^{\varepsilon}}(X^{1,\, \varepsilon}_r,\,X^{2,\,\varepsilon}_r,Y^\eps_r)\big]dr \notag \\
&+\int_0^s
[D_{x_2}V^{\varepsilon}(X^{1,\,\varepsilon}_r,X^{2,\,
\varepsilon}_r,Y^\eps_r)b^{(1)}(X^{1,\,\varepsilon}_r,\,X^{2,\,
\varepsilon}_r) - D_{y} V^\eps(X^{1,\eps}_r,X^{2,\eps}_r,Y^\eps_r)f(\frac{X^{1,\eps}_r}{\eps},
X^{2,\eps}_r,Y^\eps_r)]dr \nonumber \\
&+\int_0^s[D_{x} V^{\varepsilon}(X^{1,\,\varepsilon}_r,\,X^{2,\,\varepsilon}_r,Y^\eps_r)
\sigma (X^{1,\,\varepsilon}_r,\,X^{2,\,\varepsilon}_r)
 +
D_{y} V^\eps(X^{1,\eps}_r,X^{2,\eps}_r,Y^\eps_r)Z^\eps_r
\sigma(\frac{X^{1,\eps}_r}{\eps},X^{2,\eps}_r)]
dB_r \nonumber \\
&+\frac{1}{2}\int_0^s D_{y}^2 V^\eps(X^{1,\eps}_r,X^{2,\eps}_r,Y^\eps_r)Z^\eps_r
\sigma\sigma^\ast(\frac{X^{1,\eps}_r}{\eps},X^{2,\eps}_r)(Z^\eps_r)^\ast dr\nonumber \\
&+\frac{1}{2}\int_0^s D_{x}D_{y} V^\eps(X^{1,\eps}_r,X^{2,\eps}_r,Y^\eps_r)
\sigma\sigma^\ast(\frac{X^{1,\eps}_r}{\eps},X^{2,\eps}_r)(Z^\eps_r)^\ast dr
\end{align*}
converges to zero in probability. Let us give an explanation concerning the one but last term, which is the
most delicate one.
\begin{align*}
&\left|\int_0^sD_{y}^2 V^\eps(X^{1,\eps}_r,X^{2,\eps}_r,Y^\eps_r)Z^\eps_r
\sigma\sigma^\ast(\frac{X^{1,\eps}_r}{\eps},X^{2,\eps}_r)(Z^\eps_r)^\ast dr\right|\\
&\le C\sup_{0\le r\le s}\left|D_{y}^2 V^\eps(X^{1,\eps}_r,X^{2,\eps}_r,Y^\eps_r)\right|\text{Trace}\int_0^sZ^\eps_r
\sigma\sigma^\ast(\frac{X^{1,\eps}_r}{\eps},X^{2,\eps}_r)(Z^\eps_r)^\ast dr
\end{align*}
Since $\{\text{Trace}\int_0^sZ^\eps_r
\sigma\sigma^\ast(\frac{X^{1,\eps}_r}{\eps},X^{2,\eps}_r)(Z^\eps_r)^\ast dr,\ 0\le s\le t\}$ is the increasing process associated to a martingale  which is uniformly $L^p(\P)-$integrable for each $p\in \N$,
its $L^p(\P)$ norm is bounded, for all $p\ge1$. Finally  the same argument as above shows that
$$\sup_{0\le r\le s}\left|D_{y}^2 V^\eps(X^{1,\eps}_r,X^{2,\eps}_r,Y^\eps_r)\right|\to0$$
in probability, as $\eps\to0$.
\eop
\begin{lemma}
\label{R10} $\displaystyle \int_0^.
\bar{f}({X^{1,\,\varepsilon}_r},\,X^{2,\,\varepsilon}_r,\,Y^{\varepsilon}_r)dr
\stackrel{law}\Longrightarrow
\int_0^.\bar{f}(X^{1}_r,\,X^{2}_r,\,{Y}_r)dr$ on
$\mathcal{C}([0,\,t],\,\R)$ as $\varepsilon \longrightarrow 0$.
\end{lemma}
\noindent For the proof of this Lemma, we need
 the following two results.

\begin{lemma}
\label{AH2} Let $\displaystyle{
X^{1}_s:=x_1+\int_0^s\bar{\varphi}(X^{1}_r,\,X^2_r)dW_r,\,0\leq
s\leq t}$, and, assume {\it (A2-i), (B1)}.\\
 For $\eps>0$, let
 $\displaystyle{D^\eps_n:=\left\{s: s\in
[0,\,t]\,\,\slash \,\,|X^{1,\eps}_s|\le\frac{1}{n}\right\}}$.
\\ Define also
$\displaystyle{D_n:=\left\{s: s\in
[0,\,t]\,\,\slash \,\,|X^1_s|\le\frac{1}{n}\right\}}$.

\par\noindent
Then, there exists a constant $c>0$ such that for each $n\ge1$, $\eps>0$,
$$\E|D^\eps_n|\le\frac{c}{n}\qquad and \qquad \E|D_n|\le\frac{c}{n},$$
where $|.\,|$
denotes the Lebesgue measure.
\end{lemma}
\bop Consider the sequence $(\Psi_n)$ of functions defined as follows,
\begin{eqnarray*}
\Psi_n(x)=\left\{\begin{array}{ll}
-\frac{x}{n}-\frac{1}{2n^2} \quad if \quad x\leq -\frac{1}{n}\\\\
\frac{x^2}{2} \quad if \quad -\frac{1}{n}\leq x\leq \frac{1}{n}\\\\
\frac{x}{n}-\frac{1}{2n^2} \quad if \quad x\geq 1/n
\end{array}
\right.
\end{eqnarray*}
We put,
$\bar{\varphi}:=\bar{a}_{00}:=\rho(x_1,\,x_2)^{-1}$.

\vskip 0.2cm\noindent  Using It\^{o}'s formula, we get
\begin{eqnarray*}
\Psi_n(X^1_s)=\Psi_n(X^1_0)+\int_0^s\Psi_n^{'}(X^1_s)\bar\varphi(X^1_s,\,X^2_s)dW_s
+\frac{1}{2}\int_0^s\Psi^{"}_n(X^1_s)\bar\varphi^2(X^1_s,\,X^2_s)ds,\,s\in[0,t]
\end{eqnarray*}
Since $\bar{\varphi}$
 is lower
bounded by $C_1$, taking the expectation, we get
\begin{align*}
C_1\E\int_0^t1_{[-\frac{1}{n},\,\frac{1}{n}]}(X^1_s)ds&\leq
\E\int_0^t\Psi^{"}_n(X^1_s)\bar\varphi^2(X^1_s,\,X^2_s)ds\\
&=2\E\left[ \Psi_n(X^1_t)-\Psi_n(x_1)  \right]
\end{align*}
It follows that $\E(|D_n|)\leq
2C_1^{-1}\E\left[ \Psi_n(X^1_t)-\Psi_n(x_1)  \right]\le c/n$.
The same argument, applies to $D^\eps_n$, allows us to show the first estimate.
\eop

\begin{lemma}\label{AH20}
Consider a collection $\{Z^\eps,\ \eps>0\}$ of real valued random variables,
and a real valued random variable $Z$. Assume that for each $n\ge1$,
we have the decompositions
\begin{align*}
Z^\eps&=Z^{1,\eps}_n+Z^{2,\eps}_n\\
Z&=Z^1_n+Z^2_n,
\end{align*}
such that for each fixed $n\ge1$,
\begin{align*}
Z^{1,\eps}_n&\Rightarrow Z^1_n\\
\E|Z^{2,\eps}_n|&\le\frac{c}{\sqrt{n}}\\
\E|Z^2_n|&\le\frac{c}{\sqrt{n}}.
\end{align*}
Then $Z^\eps\Rightarrow Z$, as $\eps\to0$.
\end{lemma}

\noindent{\bf Proof.} The above assumptions
imply that the collection of random variables $\{Z^\eps,\ \eps>0\}$
is tight. Hence the result will follow from the fact that
$$\E\Phi(Z^\eps)\to\E\Phi(Z),\quad\text{as }\eps\to0$$
for all $\Phi\in C_b(\R)$ which is uniformly Lipschitz. Let $\Phi$ be such a function, and denote by $K$ its Lipschitz constant. Then
\begin{align*}
|\E\Phi(Z^\eps)-\E\Phi(Z)|&\le\E|\Phi(Z^\eps)-\Phi(Z^{1,\eps}_n)|+
+|\E\Phi(Z^{1,\eps}_n)-\E\Phi(Z^1_n)|+\E|\Phi(Z^1_n)-\Phi(Z)|\\
&\le|\E\Phi(Z^{1,\eps}_n)-\E\Phi(Z^1_n)|+2K\frac{c}{\sqrt{n}}.
\end{align*}
Hence
$$\limsup_{\eps\to0} |\E\Phi(Z^\eps)-\E\Phi(Z)|\le2K\frac{c}{\sqrt{n}},$$
for all $n\ge1$. The result follows. \eop

\noindent {\bf Proof of Lemma \ref{R10}.} For each $n\ge1$, define a
function $\theta_n\in C(\R,[0,1])$ such that $\theta_n(x)=0$ for
$|x|\le \frac{1}{2n}$, and $\theta_n(x)=1$ for $|x|\ge\frac{1}{n}$.
 We have
\begin{align*}
\int_0^t\bar{f}(X^{1,\eps}_s,X^{2,\eps}_s,Y^\eps_s)ds&=
\int_0^t\bar{f}(X^{1,\eps}_s,X^{2,\eps}_s,Y^\eps_s)\theta_n(X^{1,\eps}_s)ds
+\int_0^t\bar{f}(X^{1,\eps}_s,X^{2,\eps}_s,Y^\eps_s)[1-\theta_n(X^{1,\eps}_s)]ds\\
&=Z^{1,\eps}_n+Z^{2,\eps}_n\\
\int_0^t\bar{f}(X^{1}_s,X^{2}_s,Y_s)ds&=
\int_0^t\bar{f}(X^{1}_s,X^{2}_s,Y_s)\theta_n(X^{1}_s)ds
+\int_0^t\bar{f}(X^{1}_s,X^{2}_s,Y_s)[1-\theta_n(X^{1}_s)]ds\\
&=Z^{1}_n+Z^{2}_n
\end{align*}
Note that the mapping
$$(x^1,x^2,y)\longmapsto\int_0^t\bar{f}(x^1_s,x^2_s,y_s)\theta_n(x^1_s)ds$$
is continuous from $C([0,t])\times D([0,t])$ equipped with the product of the sup--norm and the ${\mathbf S}$ topologies into $\R$. Hence from Proposition  \ref{R8},
$Z^{1,\eps}_n\Longrightarrow Z^1_n$ as $\eps\to0$, for each fixed $n\ge1$.
Moreover, from Lemma \ref{AH2}, the linear growth property of $\bar{f}$, Lemma \ref{estX} and Proposition \ref{R5}, we deduce that
$$E|Z^{2,\eps}_n|\le\frac{c}{\sqrt{n}},\quad E|Z^{2}_n|\le\frac{c}{\sqrt{n}}.$$
 Lemma \ref{R10} now follows from Lemma \ref{AH20}.
$\blacksquare$

\vskip 0,30cm
\noindent
{\bf Proof of Proposition \ref{R11}}  Passing to the limit in the backward component of the
equation (\ref{E2}) and using Lemmas \ref{R9} and \ref{R10}, we derive assertion $(i)$.\\
Assertion $(ii)$ can be proved by using the same arguments as those in section 6
of \cite{P}. \eop

\subsection{Identification of the limit martingale. }
Since $\bar{f}$ is uniformly Lipschitz in $y$ and $H$ is bounded, then standard arguments of BSDEs (see e. g. \cite{EP2}) show that the BSDE (\ref{fbsdebar}) has a strongly unique solution and we have,
\begin{proposition}
\label{TER1} Let $(\bar Y_s,\bar Z_s,\,0\leq\,s\leq\,t)$ be the
unique solution to
BSDE (\ref{fbsdebar}). Then, for every  $s\in [0,\,t]$,\\
$\displaystyle \E|Y_s -
\bar{Y}_s|^2+\E\left([{M}-\int_0^{.}\bar{Z_r}dM^{X}_r]_t
-[{M}-\int_0^{.}\bar Z_r dM^{X}_r]_s\right)=0. $
\end{proposition}
\bop For every  $s\in [0,\,t]\setminus\mathsf{D}$,
we have
\begin{eqnarray*}
\left\{\begin{array}{l} {Y}_s=H(X_t)+\int_s^t
\bar{f}(X_r,{Y}_r)dr-({M}_t-{M}_s)\\\\
\bar Y_s = H(X_t)+\int_s^t \bar{f}(X_r,\bar Y_r)dr-\int_s^t
\bar Z_r dM^{X}_r
\end{array}
\right.
\end{eqnarray*}
Arguing as in \cite{P}, we show that $\bar M:= \int_s^.
\bar Z_r dM^{X}_r$ is a $\mathcal{F}_s$-martingale. \\ Since
$\bar f$ satisfies condition $(C1)$, we get by It\^{o}'s
formula, that
\begin{eqnarray*}
 \E|Y_s - \bar{Y}_s|^2 +
\E\left([{M}-\int_0^{.}\bar Z_r dM^{X}_r]_t
-[{M}-\int_0^{.}\bar Z_r dM^{X}_r]_s\right) \leq
C\E\int_s^t|Y_r-\bar{Y}_r|^2dr.
\end{eqnarray*}
Therefore, Gronwall's lemma yields that
$\E|Y_s-\bar{Y}_s|^2=0,\,\forall s\in[0,\,t]-\mathsf{D}$.
\\
Since
$\bar Y$ is continuous, $Y$ is c\`{a}d-lag and $\mathsf{D}$ is
countable, then $Y_s = \bar{Y}_s,\,\P$-$a.s,\,\forall
s\in[0,\,t]$. \\ Moreover, we deduce that, $\displaystyle
\E\left([{M}-\int_0^{.}\bar Z_r dM^{X}_r]_t
-[{M}-\int_0^{.}\bar Z_r dM^{X}_r]_s\right)=0. $ \eop

\vskip 0.2cm\noindent As a
consequence of Proposition \ref{TER1}, we have
\begin{corollary}
$\displaystyle
\left(Y^{\varepsilon},\,\int_0^.Z^{\varepsilon}_rdM^{X^{\varepsilon}}_r\right)
\stackrel{law}\Longrightarrow
\left(\bar Y,\,\int_0^\cdot \bar Z_r dM_r^X\right)$.
\end{corollary}

 \noindent Theorem \ref{th1} is proved.


\section{Proof of Theorem \ref{th2}.}

 Since the SDE (\ref{EH1}) is weakly unique (\cite {K}), the martingale
problem associated to $X=(X^1,\,X^2)$ is well posed.
  We then have
the following:

\begin{proposition}
\label{R4}
\ $(i)$ \ For any $t>0$, $x\in\R^d$, the
BSDE
\begin{eqnarray*}
Y^{t,\,x}_s=H(X^{x}_t)+\int_s^t
\bar f(X^{x}_r,\,Y^{t,\,x}_r)dr-\int_s^t
Z^{t,\,x}_rdM^{X^{x}}_r, 0\leq s\leq t.
\end{eqnarray*}
admits a unique solution $(Y^{t,\,x}_s,\,Z^{t,\,x}_s)_{0\leq s\leq
t}$ such that the component $(Y^{t,\,x}_s)_{0\leq
s\leq t}$ is bounded and \ $Y^{t,\,x}_0$ is deterministic.

\vskip 0.2cm \noindent $(ii)$ \ If moreover, the deterministic function,  $(t,\,x)\in
[0,\,T]\times\R^{d+1}\longmapsto v(t,\,x):=Y^{t,\,x}_0$ belongs to $ \mathcal{C}\left([0,\,T]\times
\R^{d+1},\,\R\right)$, then it  is a $L^p$-viscosity solution of the
PDE (\ref{EH2}).
\end{proposition}
\noindent
{\bf Remark.} The continuity of the map $(t,\,x)\longmapsto v(t,\,x):=Y^{t,\,x}_0$, which is assumed in assertion (ii) of Propostion \ref{R4}, will be established in Proposition \ref{R16} below.
\vskip 0.2cm\noindent
{\bf Proof of Proposition \ref{R4}.}
$(i)$ Thanks to Remark 3.5 of \cite{EP2}, it is enough to prove existence
and uniqueness for the BSDE
$$  Y^{t,\,x}_s=H(X^{x}_t)+\int_s^t
\bar f(X^{x}_r,\,Y^{t,\,x}_r)dr-\int_s^t
Z^{t,\,x}_rdB_r, 0\leq s\leq t.
$$
 Since $f$ satisfies {\bf(C)} and $\rho$ is
bounded, one can easily verify that
$\bar{f}$ is uniformly Lipschitz in $y$ uniformly with respect to $(x_1,x_2)$ and satisfies ${\bf(C1)}$-$(ii)$.
Existence and uniqueness of solution follow then from standard results for BSDEs, see e. g. \cite{EP1}.  Moreover, since $H$ is uniformly bounded and $\bar f$ satisfies the linear growth condition $\bf(C1)$-$(ii)$, one can  prove that the solution $Y^{t,\,x}$ is bounded, see e. g. \cite {B}. Finally, since $(Y_s^{t,\,x})$ is $\mathcal{F}_s^X-$adapted then $Y^{t,\,x}_0$ is measurable with respect to a trivial $\sigma-$algebra and hence it is deterministic.
\par
$(ii)$ Assume that the function
$v(t,\,x):=Y^{t,x}_0$ belongs to $ \mathcal{C}\left([0,\,T]\times
\R^{d+1},\,\R\right)$. We only prove that $v$ is a  $L^p$--viscosity sub--solution. The
proof of the super--solution property can be done similarly. Since the coefficient of  PDE under consideration are time homogeneous, then $v(t,x)$ is solution to the initial value problem (\ref{E3}) if and only if the function $u(t,x) := v(T-t,x)$ is solution to the terminal value problem.

\begin{equation}\label{semilinterminalebar}
\left\{
\begin{aligned} \frac{\partial{ u}}{\partial
t}(t,\,x)&=(\bar{L}u)(t,x)
+ \bar f(x,\,u(t,\,x)) \,\quad  t\in[0,\ T],
\\
u(T,x)&= H(x).
\end{aligned}
\right.
\end{equation}
\noindent Working with this backward PDE will simplify the details of the proofs below.

Let $X_s^{t,x}$ be the unique weak solution to SDE (\ref{EH1}). We will establish that the solution $Y$ of the Markovian BSDE
\begin{equation}
Y^{t,\,x}_s=H(X^{t,x}_T)+\int_s^T
\bar f(X^{t,x}_r,\,Y^{t,\,x}_r)dr-\int_s^T
Z^{t,\,x}_rdM_r^{X^{t,x}}, \qquad 0\leq t\leq s\leq T.
\end{equation}
define a $L^p-$viscosity sub--solution to the problem (\ref{semilinterminalebar}) by puting $u(t,x):=Y_t^{t,x}$.
\\
Let $\varphi
\in W^{1,\,2}_{p,\,loc}\left([0,\,T]\times \R^{d+1},\,\R\right)$,
let $(\widehat{t},\,\widehat{x})\in [0,\,T]\times \R^{d+1}$ be a
point which is a local maximum of  $u-\varphi$. Since $p>d+2$, then
$\varphi$ has a continuous version which we consider from now on. We
assume without loss of generality that
\begin{equation}\label{v=phi}
v(\widehat{t},\,\widehat{x})=\varphi (\widehat{t},\,\widehat{x})
\end{equation}
We will argue by contradiction.  Assume that there exists
$\varepsilon,\,\alpha>0$ such that
\begin{equation}\label{inegviscosity}
 \frac{\partial{\varphi}}{\partial
s}(s,\,x) + \bar{L}\varphi(s,x) +
\bar f(x,\,u(s,\,x))<-\varepsilon, \,\,\ \lambda\mbox{--}a. e. \mbox{ in }B_{\alpha}(\widehat{t},\,\widehat{x}).
\end{equation}
where $\lambda$ denote the Lebesgue measure.

 \noindent Since $(\widehat{t},\,\widehat{x})$ is a local maximum of $u-\varphi$, we can find a positive number $\alpha'$ (which we can  suppose equal to $\alpha$) such that
\begin{equation}
 u(t,x) \leq \varphi(t,x) \qquad \hbox{in}\quad  B_{\alpha}(\widehat{t},\,\widehat{x})
 \end{equation}
Define \begin{equation*}
\tau=\inf\left\{s\geq \widehat{t},\,;\quad
\vert X^{\widehat{t},\,\widehat{x}}_s-\widehat{x}\vert > \alpha\right\}\wedge(\widehat{t}+\alpha)
\end{equation*}
 Since $X$ is a Markov diffusion  and  $\bar f$ is uniformly Lipschitz in $y$ and satisfies condition $\bf(C1)$-$(ii)$, then arguing as in \cite{EK}, one can show that for every $r\in[\widehat{t},\ \widehat{t}+\alpha]$, \ $Y^{\widehat{t},\,\widehat{x}}_r = u(r, X^{\widehat{t},\,\widehat{x}}_r)$.  Hence, the process \  $(\bar{Y}_s,\,\bar{Z}_s):=((Y^{\widehat{t},\,\widehat{x}}_{s\wedge
\tau}),\,\1_{[0,\,\tau]}(s)(Z_s^{\widehat{t},\,\widehat{x}}))_{s\in[\widehat{t},\ \widehat{t}+\alpha]}$ \
solves the BSDE
\begin{eqnarray*}
\bar{Y}_s&=&u(\tau,\,X^{\widehat{t},\,\widehat{x}}_{\tau})
+\int_s^{\widehat{t}+\alpha}\1_{[0,\,\tau]}\,  \bar f(r,\,X^{\widehat{t},\,\widehat{x}}_r,\,u(r,\,X^{\widehat{t},\,\widehat{x}}_r)
)dr\\\\
&\qquad -&\int_s^{\widehat{t}+\alpha}\bar{Z}_r dM^{X^{\widehat{t},\,\widehat{x}}}_r,
\qquad s\in[\widehat{t},\ \widehat{t}+\alpha].
\end{eqnarray*}
On other hand, by
It\^{o}-Krylov formula, the process
$(\widehat{Y}_s,\,\widehat{Z}_s)_{s\in[\widehat{t},\ \widehat{t}+\alpha]}$, defined by $\displaystyle(\widehat{Y}_s,\,\widehat{Z}_s):=\left(\varphi(s\wedge\tau,
\,X^{\widehat{t},\,\widehat{x}}_{s\wedge\tau}),\,
\,\1_{[0,\,\tau]}(s)\nabla\varphi(s,\,X^{\widehat{t},\,
\widehat{x}}_s)\right)$ \
solves the BSDE
\begin{eqnarray*}
\widehat{Y}_s&=&\varphi(\tau,\,X^{\widehat{t},\,\widehat{x}}_{\tau})
-\int_s^{\widehat{t}+\alpha}\1_{[0,\,\tau]} [
(\frac{\partial{\varphi}}{\partial{r}}
+\bar{L}\varphi)(r,\,X^{\widehat{t},\,\widehat{x}}_r)]dr\\\\
&\qquad -&\int_s^{\widehat{t}+\alpha}\widehat{Z}_r dM^{X^{\widehat{t},\,\widehat{x}}}_r.
\end{eqnarray*}
{}From the choice of $\tau$, $(\tau,\,X^{\widehat{t},\,\widehat{x}}_{\tau})\in  B_{\alpha}(\widehat{t},\,\widehat{x})$.
Therefore, $u(\tau,\,X^{\widehat{t},\,\widehat{x}}_{\tau})\leq \varphi(\tau,\,X^{\widehat{t},\,\widehat{x}}_{\tau})$.

\noindent
Let $A := \{(t,x)\in B_\alpha(\widehat{t},\,\widehat{x}) ,\,[\frac{\partial{\varphi}}{\partial
s}+\bar{L}\varphi +
\bar f(.,u(.))](t,x)<-\varepsilon\}$ and $\bar A:= B_\alpha(\widehat{t},\,\widehat{x})\setminus A$ the complement of $A$. By (\ref{inegviscosity}), \ $\lambda (\bar A)$ = $0$.

\noindent Since the diffusion $\{X_s^{\hat t, \hat x}, s\geq t\}$ is nondegenerate, Krylov's inequality (\cite{K1}, Ch. 2, Sec. 2 $\&$ 3) implies that $\1_{\bar A}(r,\,X^{\widehat{t},\,
\widehat{x}}_r) = 0$ \ $dr\times d\P-$ $a.e.$
 It follows that
  \begin{align}\label{inegcomparaison}
 \E \int_{\widehat{t}}^{\widehat{t}+\alpha}-\1_{[0,\,\tau]} [
(\frac{\partial{\varphi}}{\partial{r}}
+\bar{L}\varphi)(r,\,X^{\widehat{t},\,\widehat{x}}_r) + \bar f(r,\,X^{\widehat{t},\,\widehat{x}}_r,\,u(r,\,X^{\widehat{t},\,
\widehat{x}}_r)
)]
) dr
 \geq \ \E(\tau-\widehat{t})\, \varepsilon
 > \ 0
\end{align}
This implies that \ $[- \1_{[0,\,\tau]}  [
(\frac{\partial{\varphi}}{\partial{r}}
+\bar{L}\varphi)(r,X^{\widehat{t},\widehat{x}}_r) + \bar f(r,X^{\widehat{t},\widehat{x}}_r,u(r,\,X^{\widehat{t},
\widehat{x}}_r)
)]
)] \ > \,0$ \ on a set of $dt\times d\P$ positive measure. Therefore,
 the strict comparison theorem (Remark 2.5 in
\cite{EP2}) shows that
$\bar{Y}_{\widehat{t}}<\widehat{Y}_{\widehat{t}}$ , that is
$u(\widehat{t},\,\widehat{x})<\varphi(\widehat{t},\,\widehat{x})$,
which contradicts our assumption  (\ref{v=phi}).


\vskip 0.3cm Under assumptions {\bf (A), (B)}, the SDE $(\ref{EH1})$ has a unique weak solution, see \cite{K}. We then have the following continuity property.
\begin{proposition}(Continuity in law of the map $x\mapsto X_.^x$)
\vskip 0,2cm
\label{AH1}
\noindent Assume {\bf (A), (B)}.
Let $X_s^x$ be the unique weak solution of the SDE $(\ref{EH1})$,
and
\\
 $\displaystyle
X^{n}_s:=x_n+\int_0^s
\bar{b}(X^{n}_r)dr+\int_0^s\bar{\sigma}(X^{n}_r)dB_r,\,\,0\leq
s\leq t$\\
Assume that $x_n\to x=(x^1,\,x^2)\in
\R^{1+d}$ as $n\to\infty$.
Then
\noindent  $\displaystyle X^{n}\stackrel{law}\Longrightarrow X^x$.
\end{proposition}
\bop  Since $\bar{b}$ and $\bar{\sigma}$ satisfy
{\bf (A), (B)}, one can easily check that the sequence $X^n$ is tight in $\mathcal{C}([0,\,t]\times\R^{d+1})$.
By Prokhorov's theorem, there exists a subsequence (denoted also by $X^n$) which converges weakly
 to a process $\widehat X$.
 We shall show that $\widehat X$ is a weak solution of SDE $(\ref{EH1})$.\\
$\bullet$ {\it Step 1:} For every $\varphi\in C_c^{\infty}(\R^{1+d})$,
\begin{eqnarray*}
\forall
u\in[0,\,t],\quad\varphi(\widehat{X}_r)-\int_0^u\bar{L}\varphi(\widehat{X}_v)dv\quad
\mbox{ is a $\mathcal{F}^{\widehat{X}}$-martingale.}
\end{eqnarray*}
 All we need to show is that for every  $\varphi\in C_c^{\infty}(\R^{1+d})$,  every $0\leq s\leq u$ and every
function $\Phi_s$ of $(X^{x_n}_r)_{0\leq r<s}$ which is bounded and
continuous for the topology of uniform convergence, as $n\to\infty$,
\begin{eqnarray*}
0=\E\left\{[\varphi(X^{x_n}_r)-\varphi(X^{x_n}_s)-\int_s^r\bar{L}\varphi(X^{x_n}_\alpha)d\alpha]\Phi_s(X^{x_n}_.)\right\}\\\\
\longrightarrow
\E\left\{[\varphi(\widehat{X}_r)-\varphi(\widehat{X}_s)-\int_s^r\bar{L}\varphi(\widehat{X}_\alpha)d\alpha]\Phi_s(\widehat{X}_.)\right\}
\end{eqnarray*}
Indeed, since $\varphi,\,\Phi$ are continuous functions and $\bar{L}\varphi$
is continuous away from the set $\{x_1=0\}$, similar argument as that developed in the proof of Lemma \ref{R10} gives
\begin{eqnarray*}
[\varphi(X^{x_n}_r)-\varphi(X^{x_n}_s)-\int_s^u\bar{L}\varphi(X^{x_n}_v)dv]
\Phi_s(X^{x_n}_.) \stackrel{law}\longrightarrow
[\varphi(\widehat{X}_r)-\varphi(\widehat{X}_s)-\int_s^u\bar{L}\varphi(\widehat{X}_v)dv]\Phi_s(\widehat{X}_.)
\end{eqnarray*}
Since $\varphi,\,\Phi$ are bounded functions, $\bar{L}\varphi$ has at most linear growth at infinity and
$$\sup_n\E(\sup_{s\in[0,\,t]}|X^{x_n}|^2)<\infty,$$ the result
follows by uniform integrability. Hence
$$
\E\left\{[\varphi(\widehat{X}_r)-\varphi(\widehat{X}_s)
-\int_s^u\bar{L}\varphi(\widehat{X}_v)dv]\Phi_s(\widehat{X}_.)\right\}=0$$
and therefore
$\varphi(\widehat{X}_r)-\varphi(\widehat{X}_s)-\int_s^r\bar{L}\varphi(\widehat{X}_v)dv$
is a $\mathcal{F}^{\widehat{X}}_r$--martingale.

\vskip 0,15cm
\noindent $\bullet${\it Step 2:} From {\it step 1}, there
exists a $\mathcal{F}^{\widehat{X}}$-Brownian motion $\widehat{B}$
such that,
\begin{eqnarray*}
\widehat{X}_s=x+\int_0^s
\bar{b}(\widehat{X}_r)dr+\int_0^s\bar{\sigma}(\widehat{X}_r)d\widehat{B}_r, \quad 0\leq s\leq t.
\end{eqnarray*}
Weak uniqueness of the SDE (\ref{EH1}) allows us to deduce that
$\widehat{X}= X^x$ in law sense.
\eop
\begin{proposition}
\label{R16}
 Assume {\bf (A), (B), (C)}. Then,
\vskip 0,2cm \noindent $\displaystyle(i)\, \lim_{\varepsilon
\rightarrow 0}Y^{\varepsilon}_0$ = $Y_0^{t,x}$.

\vskip 0.2cm \noindent $\displaystyle(ii)$ The map $(t,x)\longmapsto
  Y^{t,\,x}_0$ is continuous.

\vskip 0.2cm \noindent $\displaystyle(iii)$ For $p>d+2$, the function $v(t,x) := Y^{t,\,x}_0$ is
a $L^{p}$-viscosity solution to the PDE  (\ref{E3}).
\end{proposition}
\bop
 $(i)$ Let $Y^{t,x}$ be the limit process defined in Proposition \ref{R8}.
 We have
\begin{eqnarray*}
\left\{\begin{array}{l}
Y^{\varepsilon}_0=H(X^{\varepsilon}_t)+\displaystyle\int_0^t f(\frac{{X}^{1,\varepsilon}_r}{\varepsilon} \ , X^{2,\,\varepsilon}_r,\,Y^{\varepsilon}_r)dr-M^{\varepsilon}_t \\\\
{Y_0^{t,x}}=H(X^x_t)+\displaystyle\int_0^t\bar{f}(X^x_r, {Y_r^{t,x}})dr-{M}_t
\end{array}
\right.
\end{eqnarray*}
{}From Jakubowski \cite{J}, the projection: $y\mapsto y_t$ is
continuous from $\mathcal{D}([0,\,t];\R)$ into $\R$ for the $\bf{S}$--topology. We then deduce from the convergence of the above right--hand sides that
$Y^{\varepsilon}_0$ converges towards ${Y}_0$ in distribution. Since
$Y^{\varepsilon}_0$ and $Y_0$ are deterministic, this means exactly that
$Y^{\varepsilon}_0\to Y_0$

\vskip 0,3cm
\noindent
$(ii)$
Let  $(t_n,\,x_n)\rightarrow (t,\,x)$. We assume that $t>t_n>0$. We have,
\begin{eqnarray}
\label{FIN1}
Y^{t_n,\,x_n}_s&=&H(X^{x_n}_{t_n})+\int_s^{t_n}
\bar{f}(X^{x_n}_r,\,Y^{t_n,\,x_n}_r)dr-\int_s^{t_n}
Z^{t_n,\,x_n}_rdM^{X^{x_n}}_r,\,\,0\leq\,s\leq\,t_n,
\end{eqnarray}
where $X^{x_n}\stackrel{law}\Rightarrow\,X^x$. \\
Since $H$ is a bounded continuous function
and $\bar{f}$ satisfies $(C1)$, one can easily
show that the sequence $\{(Y^{t_n,\,x_n},\,\int_0^.1_{[s,t_n]}(u)Z^{x_n}_rdM^{X^{x_n}}_r)\}_{n\in\N^*}$ is tight in $\mathcal{D}([0,\,t];\R^2)$.\\
Let us rewrite the equation (\ref{FIN1}) as follows
\begin{eqnarray}
\label{FIN2}
Y^{t_n,\,x_n}_s&=&H(X^{x_n}_{t_n})+\int_s^{t}
\bar{f}(X^{x_n}_r,\,Y^{t_n,\,x_n}_r)dr-\int_s^{ t}1_{[s,t_n]}(u)
Z^{t_n,\,x_n}_rdM^{X^{x_n}}_r\\\nonumber
&-&\int_{t_n}^t\bar{f}(X^{x_n}_r,\,Y^{t_n,\,x_n}_r)dr,\,0\leq\,s\leq t.\\\nonumber
&=&A^1_n+A^2_n\nonumber
\end{eqnarray}
\vskip 0,15cm
\noindent $\bullet$ {\it Convergence of $A^2_n$}\\
Since $\bar{f}$ is bounded, $\displaystyle
\E\left|\int_{t_n}^t\bar{f}(X^{x_n}_r,\,Y^{t_n,\,x_n}_r)dr\right|
\leq K |t-t_n|$. Hence $A^2_n$ tends to zero in probability.\\
\noindent $\bullet$ {\it Convergence of $A^1_n$}\\
Denote by  $(Y',\,M')$ the weak limit of $\{(Y^{t_n,\,x_n},\,\int_0^.1_{[s,t_n]}(u)Z^{x_n}_rdM^{X^{x_n}}_r)\}_{n\in\N^*}$.
The same proof as that of Lemma \ref{R10} establishes that $\displaystyle \int_s^{t}
\bar{f}(X^{x_n}_r,\,Y^{t_n,\,x_n}_r)dr\stackrel{law}\Longrightarrow \int_s^{t}
\bar{f}(X^{x}_r,\,Y'_r)dr$.\\
\vskip 0,15cm
\noindent Passing to the limit in (\ref{FIN2}), we obtain that
\begin{eqnarray*}
Y'_s&=&H(X^{x}_{t})+\int_s^{t}
\bar{f}(X^{x}_r,\,Y'_r)dr-(M'_{t}-M'_s),\,\,s\in[0,\,t]\cap D^c.
\end{eqnarray*}
The uniqueness of the considered BSDE ensures that $\forall s\in[0,\,t],\,Y'_s=Y_s^{t,\,x}\,\P$-ps.
Hence $\displaystyle Y^{t_n,\,x_n}\stackrel{law}\Rightarrow\,Y^{t,\,x}$.
As in $(i)$, one derive that $Y^{t_n,\,x_n}_0\stackrel{law}\Longrightarrow\,Y^{t,\,x}_0$
which yields to the continuity
of $Y_0^{t,\,x}$.

\noindent Assertion $(iii)$ follows from $(ii)$ and the second statement of Proposition \ref{R4}. \eop



\appendix \section{Appendix: S-topology}
The $\bf{S}$--topology has been introduced by Jakubowski (\cite{J},
1997) as a topology defined on the Skorohod space of c\`adl\`ag
functions: $\mathcal{D}([0,\,T];\,\R)$. This topology is weaker
than the Skorohod topology but tightness criteria are easier to
establish. These criteria are the same as the one used in
Meyer-Zheng \cite{MZ}. \\Let $N^{a,\,b}(z)$
denotes the number of up-crossing of the function
$z\in\mathcal{D}([0,\,T];\,\R)$ from level $a$ to level $b$ ($a<b$). We recall
some facts about the $\bf{S}$--topology.
\begin{proposition}(A criteria for S-tight). A sequence $(Y^{\varepsilon})_{\varepsilon>0}$ is said to be
$\bf{S}$--tight if and only if it is relatively compact for the $\bf{S}$--topology.\\
Let $(Y^{\varepsilon})_{\varepsilon>0}$ be a family of stochastic
processes in $\mathcal{D}([0,\,T];\,\R)$. Then this family is
tight for the $\bf{S}$--topology if and only if
$(\|Y^{\varepsilon}\|_{\infty})_{\varepsilon>0}$ and
$(N^{a,\,b}(Y^{\varepsilon}))_{\varepsilon>0}$ are tight for each
$a<b$.
\end{proposition}
\noindent Let
$\left(\Omega,\,\mathcal{F},\,\P,\,(\mathcal{F}_t)_{t\geq
0}\right)$ be a stochastic basis. If $(Y)_{0\leq t\leq T}$ is a
process in $\mathcal{D}([0,\,T];\,\R)$ such that $Y_t$ is
integrable for any $t$, {\it the  conditional variation of $Y$} is
defined by
$$
CV(Y)=\sup_{n\ge1,\ 0\leq t_1<...<t_n=T}
\sum_{i=1}^{n-1}\E[|\E[Y_{t_{i+1}}-Y_{t_i}\left|\right.\mathcal{F}_{t_i}]|].
$$
The process $Y$ is called a $quasimartingale$ if $CV(Y)<+\infty$. When $Y$
is a $\mathcal{F}_t$-martingale, $ CV(Y)=0$. A variation of Doob's
inequality (cf. lemma 3, p. 359 in Meyer and Zheng \cite{MZ}, where
it is assumed that $Y_T=0$) implies that
$$
\P\left[\sup_{t\in[0,\,T]}|Y_t|\geq k\right]\leq\frac{2}{k}\left
(CV(Y)+\E\left[\sup_{t\in[0,\,T]}|Y_t|\right] \right),
$$
$$
\E\left[N^{a,\,b}(Y)\right]\leq\frac{1}{b-a}\left
(|a|+CV(Y)+\E\left[\sup_{t\in[0,\,T]}|Y_t|\right]\right ).
$$
It follows that a sequence $(Y^{\varepsilon})_{\varepsilon>0}$ is
$\bf{S}$-tight whenever
$$
\sup_{\varepsilon>0}\left(CV(Y^{\varepsilon})+\E\left[\sup_{t\in[0,\,T]}|Y_t^{\varepsilon}|\right]\right)<+\infty.
$$
\begin{theorem}
\label{B1} Let $(Y^{\varepsilon})_{\varepsilon>0}$ be a $\bf{S}$-tight
family of stochastic process whose trajectories belong to $\mathcal{D}([0,\,T];\,\R)$. Then
there exists a sequence $(\varepsilon_k)_{k\in\N}$ decreasing to
zero, some process $Y\in \mathcal{D}([0,\,T];\,\R)$ and a
countable subset $D\in [0,\,T]$ such that for any $n\ge1$ and any
$(t_1,\,...,\,t_n)\in[0,\,T]\backslash D$,
$$
(Y^{\varepsilon_k}_{t_1},\,...,\,Y^{\varepsilon_k}_{t_n})\stackrel{\mathcal{D}ist}{\longrightarrow}
(Y_{t_1},\,...,\,Y_{t_n})
$$
\end{theorem}
\begin{remark}
\label{B2} The projection
$\pi_T\, : \,y\in(\mathcal{D}([0,\,T];\,\R),\,{\bf S})\mapsto\,y(T)$ is
 continuous (see Remark 2.4, p.8 in Jakubowski \cite{J}), but $y\mapsto\,y(t)$ is not continuous for each $0\leq t\leq T$.
\end{remark}
\begin{lemma}
\label{B3} Let $(U^{\varepsilon},\,M^{\varepsilon})$ be a
multidimensional process in
$\mathcal{D}([0,\,T];\,\R^{p})\,(p\in\N^{*})$ converging to
$(U,\,M)$ in the S-topology. Let
$(\mathcal{F}_t^{U^{\varepsilon}})_{t\geq 0}$ (resp.
$(\mathcal{F}_t^{U})_{t\geq 0}$) be the minimal complete admissible
filtration generated by $U^{\varepsilon}$ (resp. $U$). We assume
moreover that for every $T>0$,
$\sup_{\varepsilon>0}\E\left[\sup_{0\leq t\leq
T}|M^{\varepsilon}_t|^2\right]<C_T$.  \\
 If  $M^{\varepsilon}$
 is a $\mathcal{F}^{U^{\varepsilon}}$-martingale and $M$ is
$\mathcal{F}^{U}$-adapted, then $M$ is a
$\mathcal{F}^{U}$-martingale.
\end{lemma}
\begin{lemma}
\label{B4} Let $(Y^{\varepsilon})_{\varepsilon>0}$ be a sequence
of process converging weakly in $\mathcal{D}([0,\,T];\,\R^p)$  to
$Y$. We assume that $\sup_{\varepsilon>0}\E\left[\sup_{0\leq t\leq
T}|Y^{\varepsilon}_t|^2\right]<+\infty$. Then for any $t\geq
0,$ $E\left[\sup_{0\leq t\leq T}|Y_t|^2\right]<+\infty$.
\end{lemma}

\noindent
$\bf{Acknowlegement.}$ \
The authors are sincerely grateful to the referee for many
useful remarks which have leaded to an improvement of the paper.

The second author thanks IMATH laboratorie of universit\'e du Sud
Toulon-Var and the  LATP laboratory of universit\'e de Provence,  Marseille, France,
for their kind hospitality.

\end{document}